\documentclass[12pt]{article}

\newtheorem{definition}{Definition}

\usepackage{epsfig}
\usepackage{url}
\usepackage{amssymb}
\usepackage{amsmath}
\usepackage{amsfonts}

\def\Dbar{\leavevmode\lower.6ex\hbox to 0pt{\hskip-.23ex \accent"16\hss}D}

\def\bZ{{\mbox{\bf Z}}}

\newcommand{\nc}{\newcommand}
\nc{\cP}{{\cal P}}

\begin{document}

{\bf\LARGE
\begin{center}
Symmetric Hadamard matrices of orders 268, 412, 436 and 604
\end{center}
}

{\Large
\begin{center}
N.A. Balonin\footnote{Saint-Petersburg State University of Aerospace Instrumentation, 67, B. Morskaia St., 190000, Saint-Petersburg, Russian Federation, E-mail: korbendfs@mail.ru} and 
D.{\v{Z}}. {\Dbar}okovi{\'c}\footnote{University of Waterloo,
Department of Pure Mathematics and Institute for
Quantum Computing, Waterloo, Ontario, N2L 3G1, Canada
e-mail: \url{djokovic@uwaterloo.ca}}
\end{center}
}

\begin{abstract}
We construct many symmetric Hadamard matrices of small order 
by using the so called propus construction. The necessary 
difference families are constructed by restricting the search to the families which admit a nontrivial multiplier. Our main result is that we have constructed, for the first time, 
symmetric Hadamard matrices of order 268, 412, 436 and 604. \end{abstract}



\section{Introduction}

The construction of symmetric Hadamard matrices was stagnating for long time while that of skew-Hadamard matrices advanced 
rapidly. The reason for this discrepancy was the fact that 
for the latter we had a very versatile tool, namely the 
Goethals-Seidel (GS) array, while for the former such tool 
was missing. The new tool for the construction of the 
symmetric Hadamard matrices, so called propus array, was discovered recently \cite{SB:2017} by J. Seberry and the 
second author. It was already used in 
\cite{BBDKM:2017,BDK:2018,Djokovic:SpecMatC:2015}
to construct many propus Hadamard matrices (such matrices are always symmetric) including some having new orders. 

The authors of \cite{SB:2017} observed that 
the well known Turyn series of Williamson quadruples (of  
symmetric circulant blocks) gives the first infinite series 
of propus Hadamard matrices. They also give a variation of the 
propus array in which they plug symmetric and commuting 
Williamson type quadruples to construct another infinite series 
of symmetric Hadamard matrices. Yet another infinite series 
of propus Hadamard matrices was identified in 
\cite[Theorem 5]{Djokovic:SpecMatC:2015}.

In this paper we continue our previous work 
\cite{BBDKM:2017,BDK:2018} where we used the propus 
construction to find new symmetric Hadamard matrices. 
We refer to these papers and \cite{DK:GS:2018} for the more comprehensive description of this construction and the definitions of the GS-array and GS-difference families. 
As the propus difference families play a crucial role in the 
paper, we shall define them precisely in the next section and 
specify the propus array that we use.

The first Hadamard matrix of order $4\cdot 67=268$ was 
constructed by Sawade in 1985 \cite{Saw}. The first 
skew-Hadamard matrix of the same order was constructed 
in 1992 by one of the authors \cite{Djokovic:BAMS:1992}. 
However a symmetric Hadamard matrix of order $268$ was not 
discoverd so far. We present in Sect. \ref{sec:67} six propus difference families in the cyclic group $\bZ_{67}$ which we use to construct six symmetric Hadamard matrices of order $268$. 
Moreover, in the same section we also 
construct the first examples of symmetric Hadamard matrices 
of orders $412$, $436$ and $604$.
Examples of symmetric Hadamard matrices of order $4v$ are 
now known 
\cite{BBDKM:2017,BDK:2018,CK-Had:2007,Djokovic:SpecMatC:2015}
for all odd positive integers $v<200$ except for 
$$
59,65,81,89,93,101,107,119,127,133,149,153,163,167,179,183,189,
191,193.
$$

The binary sequences, i.e., $\{\pm1\}$-sequences, of 
length $v \equiv 1 \pmod{4}$ are called {\em optimal} if 
the off-pick values of its periodic autocorrelation 
function are $+1$ or $-3$. Such sequence is {\em balanced} 
if its sum is $\pm1$. A computer generated list of binary 
balanced optimal sequences of length $v\equiv1 \pmod{4}$ 
is given in \cite{Arasu:2011} for $v\le47$. As a byproduct of our computations of propus difference families we have obtained 
binary balanced optimal sequences of lengths 49 and 61. 
They are presented in Sect. \ref{sec:bobs}.

In addition to the propus difference families
used in Sect. \ref{sec:67}, we give a more extensive list of such families in Sect. \ref{sec:PropDF}.

While trying to verify the proof of \cite[Corollary 1]{SB:2017} we observed that this corollary is stated incorrectly. 
The second sentence of the corollary should read:
``Then there exist symmetric Williamson type matrices 
of order $q+2$ and a symmetric propus-type Hadamard matrix 
of order $4(q+2)$.'' 
Consequently, $4(2q+1)$ should be replaced with $4(q+2)$ in
the abstract as well as in line 3 on p. 351. 
Further, the two lists, one on p. 352 and the other on p. 356 
should be corrected. The integers $59,67,81,89,105,111,119,127$ 
should be removed from the former, while $97,99$ should be 
removed from and $59,67,89,119,127$ inserted into the latter. 
(The cases $59,89,119,127$ are still unresolved.)

\section{Preliminaries} \label{sec:prelim}

Let $G$ be a finite abelian group of order $v>1$. 
Let $(X_i)$, $i=1,2,\ldots,m$, be a difference family in $G$. 
We fix its parameter set 

\begin{equation} \label{param-gen}
(v;k_1,k_2,\ldots,k_m;\lambda), \quad k_i=|X_i|.
\end{equation}

Recall that these parameters satisfy the equation
\begin{equation} \label{lambda}
\sum_{i=1}^t k_i(k_i-1)=\lambda(v-1).
\end{equation}

The set of difference families in $G$ having this parameter set is invariant under the following elementary transformations:

(a) For some $i$ replace $X_i$ by a translate $g+X_i$,
$g\in G$.

(b) For some $i$ replace $X_i$ by $-X_i$.

(c) For all $i$ replace $X_i$ by its image $\alpha(X_i)$
under an automorphism $\alpha$ of $G$.

(d) Exchange $X_i$ and $X_j$ provided that $|X_i|=|X_j|$.

\begin{definition} \label{ekv}
We say that two difference families with the same parameter set are {\em equivalent} if one can be transformed to the other by a finite sequence of elementary transformations.
\end{definition}

\begin{definition} \label{mult}
Let $(X_i)$ be a difference family in $G$. We say that an 
automprphism $\alpha$ of $G$ is a {\em multiplier} of this 
family if each set $\alpha(X_i)$ is a translate of $X_i$.
\end{definition}

If a positive integer $m$ is relatively prime to $v$ then 
the multiplication by $m$ is an automorphism of $G$. If this 
automorphism is a multiplier of a difference family, then we 
also say that the integer $m$ is a {\em multiplier} or a 
{\em numeric multiplier} of that family.

The multipliers of a difference family in $G$ form a subgroup 
of the automorphism group of $G$. All difference families that 
we construct in this paper have nontrivial multipliers. This 
follows from the fact that we use only the base blocks $X_i$ 
which are union of orbits of a fixed nontrivial subgroup $H$ 
of the automorphism group of $G$. We refer to this method of 
constructing difference families as the {\em orbit method}.

We are only interested in Goethals-Seidel (GS) difference families formally introduced in 
\cite{DK:GS:2018} and \cite{DK:KonAbGr:2018}. 
They consist of four base blocks $(X_1,X_2,X_3,X_4)$ 
and their parameter sets, also known as the 
{\em  GS-parameter sets}, satisfy besides the obvious condition (\ref{lambda}) (with $m=4$) also the condition

\begin{equation} \label{sum-k}
\sum_{i=1}^4 k_i=\lambda+v.
\end{equation}

By eliminating the parameter $\lambda$ from the equations 
(\ref{lambda}) and (\ref{sum-k}), we obtain that

\begin{equation} \label{sum-kv}
\sum_{i=1}^4 (v-2k_i)^2=4v.
\end{equation}

If $k_i=k_j$ for some $i\ne j$ in a GS-parameter set 
$(v;k_1,k_2,k_3,k_4;\lambda)$ then we say that this 
parameter set is a {\em propus parameter set}.

In fact we shall use only a very special class of GS-difference 
families known as {\em propus difference families}. We adopt here the following definiton of these families.

\begin{definition} \label{df:propus}
A {\em propus difference family} is a GS-difference family 
$(X_i)$, $i=1,2,3,4$, subject to two additional conditions:

(a) two of the base blocks are equal, say $X_i=X_j$
for some $i<j$, which implies that $k_i=k_j$;

(b) at least one of the other two base blocks is symmetric.
\end{definition}

(We say that a subset $X\subseteq G$ is {\em symmetric}
if $-X=X$.)

Unless stated otherwise, we shall assume from now on that 
$G$ is cyclic. We identify $G$ with the additive group of the 
ring $\bZ_v$ of integers modulo $v$. We denote by $\bZ_v^*$ 
the group of units (invertible elements) of $\bZ_v$. 
We identify the automorphism group of $G$ with $\bZ_v^*$. 
Thus, every automorphism $\alpha$ of $\bZ_v$ is just the multiplication modulo $v$ by some integer $k$ relatively prime 
to $v$.

To any subset $X\subseteq\bZ_v$ we associate the binary sequence
(i.e., a sequence with entries $+1$ and $-1$) of length $v$, say
$(x_0,x_1,\ldots,x_{v-1})$, where $x_i=-1$ if and only
if $i\in X$. By abuse of language, we shall use the
symbol $X$ to denote also the binary sequence
associated to the subset $X$.

Let $(X_i)$ be a GS-difference family in $\bZ_v$.
Further, let $A_i$ be the circulant matrix having the sequence $X_i$ as its first row. Then the $A_i$ satisfy the equation

\begin{equation} \label{sum-A}
\sum_{i=1}^4 A_i^T A_i=4vI_v,
\end{equation}
where $I_v$ is the identity matrix of order $v$.
This equation guarantees that, after plugging the
$(A_i)$ into the GS-array, we obtain a Hadamard matrix.

If $(X_i)$ is a propus difference family, we say that the 
corresponding matrices $(A_i)$ are {\em propus matrices}. 
By plugging these $(A_i)$, in suitable order, into the 
{\em propus array}
\begin{equation} \label{Propus-array}
\left[ \begin{array}{cccc}
-A_1 & A_2R & A_3R & A_4R \\
A_3R & RA_4 & A_1 & -RA_2 \\
A_2R & A_1 & -RA_4 & RA_3 \\
A_4R & -RA_3 & RA_2 & A_1
\end{array} \right],
\end{equation}
where $R$ is the back-diagonal permutation matrix, we obtain a 
symmetric Hadamard matrix of order $4v$. The ordering should 
be  chosen so that $A_1$ is symmetric and $A_2=A_3$. 

We construct the base blocks $X_i$ as unions of certain 
orbits of a small nontrivial subgroup $H\subseteq\bZ^*_v$ 
(mostly of order 3 or 5). When recording a base block, 
to save space, we just list the representatives of the orbits 
which occur in the block. As a representative, we always choose 
the smallest integer of the orbit.

\section{The cases $v=67,103,109,151$} \label{sec:67}

In this section we list six non-equivalent examples 
of propus difference families in $\bZ_{67}$, three such 
families in $\bZ_{103}$, two in $\bZ_{109}$, and 
a single one in $\bZ_{151}$. 
By using the propus array, they provide the first examples of 
symmetric Hadamard matrices of orders 268, 412, 436 and 604, 
respectively.

In the case $v=67$, up to a permutation of the $k_i$s, there 
are three feasible propus parameter sets for the subgroup 
$H=\{1,29,37\}\subseteq\bZ^*_{67}$. For each of them we have found several propus difference families. We list only two 
families per parameter set. The block $X_4$ is symmetric in 
the first two families while $X_1$ is symmetric in the 
remaining four families. 

Let us explain how we record the base blocks. 
As an example, we take the block $X_2$ of the first family in 
Table 1. It is the union of ten $H$-orbits whose 
representatives are the integers $0,2,4,6,16,17,25,27,30,41$. 
As each nontrivial orbit has size 3, the block $X_2$ has 
the size $1+9\cdot 3=28$. 
The blocks $X_1$ and $X_4$ are given similarly. 
In all difference families listed in this and the next section 
we have $X_2=X_3$ and we record only the blocks 
$X_1$, $X_2$ and $X_4$ in that order. The families having the same parameter set are separated by a semicolon. 

For the cases $v=103$ and $v=109$ we use again the subgroups 
$H$ of order 3, namely $\{1,46,56\}\subset\bZ_{103}^*$ and 
$\{1,45,63\}\subset\bZ_{109}^*$. For $v=103$ we found two 
non-equivalent propus difference families having the 
same parameter set and for $v=109$ we found three such families. In all six families the block $X_4$ is symmetric. 

For the case $v=151$ we use the subgroup of order five. Only 
one propus difference family was found. The symmetric block is
$X_1$ .

\begin{center} Table 1. Propus difference families in 
$\bZ_{67}$, $\bZ_{103}$, $\bZ_{109}$  and $\bZ_{151}$. \end{center}

\begin{verbatim}
(67;33,28,28,31;53), H={1,29,37}
[1,3,4,10,12,15,17,30,34,36,41],[0,2,4,6,16,17,25,27,30,41],
[0,1,4,5,8,10,16,18,30,32,36];
[1,2,8,15,16,18,25,30,32,34,36],[0,2,3,6,8,9,17,18,34,36],
[0,1,2,4,5,9,16,17,18,30,41]

(67;30,31,31,27;52), H={1,29,37}
[1,5,6,15,16,17,27,30,34,41],[0,2,4,9,10,12,16,23,30,36,41],
[5,8,9,12,16,17,23,25,41];
[3,5,8,10,12,16,23,25,32,36],[0,5,6,9,12,15,16,17,23,27,30],
[1,2,3,4,8,27,30,32,36]

(67;30,30,30,28;51), H={1,29,37}
[3,4,5,8,10,16,18,23,32,36],[3,6,9,10,12,15,17,23,25,41],
[0,5,9,10,12,15,17,27,30,41];
[2,3,4,9,10,17,18,23,32,41],[1,2,9,16,17,23,27,32,34,41],
[0,3,10,15,16,17,23,27,32,34]

(103;48,51,51,42;89), H={1,46,56}
[3,4,14,17,19,21,29,30,31,33,38,40,49,51,55,62],
[2,3,4,6,7,14,15,22,29,30,31,38,42,44,47,49,62],
[3,6,8,10,15,17,21,31,33,38,42,44,55,60];
[1,3,6,8,10,11,21,30,33,40,44,47,49,51,55,62],
[5,6,7,11,12,14,19,23,29,30,38,40,47,51,55,60,62],
[4,6,7,8,10,12,17,20,22,33,42,44,49,55]

(109;52,49,49,48;89), H={1,45,63}
[0,3,4,6,9,10,11,12,18,19,20,24,31,36,43,48,50,60],
[0,1,2,3,5,9,10,16,19,20,23,25,41,46,55,57,62],
[1,2,4,6,9,10,15,19,20,24,31,36,38,46,48,57];
[0,3,5,8,11,12,13,15,18,20,30,31,41,43,46,53,55,57],
[0,1,2,3,5,8,11,12,13,16,29,31,38,41,48,50,57],
[3,6,8,10,18,20,23,24,25,29,41,48,55,57,60,62];
[0,1,2,3,6,9,10,12,15,18,24,25,36,41,43,48,53,57],
[0,1,3,6,8,9,11,12,13,18,23,29,31,36,41,43,57],
[1,3,9,11,13,16,18,29,30,31,43,46,50,53,62,67]

(151;71,71,71,66;128), H={1,8,19,59,64}
[0,2,5,6,7,11,15,17,23,27,30,34,37,51,68],
[0,1,2,3,4,14,17,23,27,28,34,47,51,68,87],
[0,1,2,3,4,5,7,10,29,34,46,47,51,68]
\end{verbatim}

\section{Some new balanced optimal binary sequences}
\label{sec:bobs}

In this section we list some balanced optimal binary 
sequences of lengths 49 and 61. They arose as a byproduct 
of our search for propus difference families. We say that 
a binary sequence of length $v$ has {\em three-level 
autocorrelation function} if this function takes exactly 
three distinct values, including the value $v$ at shift 0.

Up to a permutation of the $k_i$s, there are three feasible propus parameter sets for the 
subgroup $H=\{1,18,30\}$ of $\bZ^*_{49}$. 
We discard the one with all $k_i=21$ as it probably does not 
admit any propus difference family, see \cite{BDK:2018}.
In Table 2 we list five propus difference families for $v=49$ 
and a single family for $v=61$.

\begin{center} Table 2. Three-level autocorrelation functions from propus difference families. 
\end{center}

\begin{verbatim}
(49;22,24,24,18;39), H={1,18,30}
[0,1,6,7,8,9,13,16],[3,7,8,9,13,16,21,29],[3,6,8,12,16,29];
[0,2,7,8,13,16,19,26],[2,6,9,12,16,24,26,29],[1,3,7,8,19,21];
[0,1,3,4,12,13,16,24],[1,6,8,13,16,19,24,29],[1,4,6,16,19,26]

(49;22,22,22,19;36), H={1,18,30}
[0,4,6,7,9,13,19,26],[0,1,6,7,9,12,16,29],[0,1,6,7,16,19,21];
[0,3,4,6,7,12,19,29],[0,1,2,4,7,8,13,19],[0,1,3,7,8,19,21]

(61;25,30,30,25;49), H={1,13,47}
[0,6,8,11,16,18,23,32,36],[1,2,3,9,12,22,27,28,31,36],
[0,4,7,8,9,11,16,27,28]
\end{verbatim}

The block $X_2$, of cardinality 24, in the first example is 
$$ 
X_2=\{3,5,7,8,9,13,14,15,16,21,25,28,29,32,35,37,38,39,41,42,
43,44,46,47\}.
$$
The values of the periodic autocorrelation function of the 
corresponding sequence $X_2$, for the shifts in the range 
$0,1,\ldots,24$, are:
$$
49,1,-3,-3,1,-3,1,1,-3,-3,1,-3,-3,-3,1,-3,1,-3,1,1,-3,1,1,1,-3.
$$
Thus the correlation values of $X_2$ occupy just three levels 
49, 1 and $-3$. In the terminology of \cite[p. 144]{Arasu:2011} 
(see also \cite{MB:1998}) 
the sequence $X_2$ is a balanced optimal binary sequence of length $49$. Such sequences of lengths $v\equiv 1 \pmod{4}$ 
are listed there on the same page for $v\le 45$. Our 
sequence $X_2$ extends that list one step further. 
The meaning of the word `balanced' in this context 
is that the sum of the sequence is 1 or $-1$.

The sequences $X_2$ in the second and third example also 
have only 3 correlation values but this time these values 
are 49, 1 and $-7$ and so they are not optimal.

The block $X_2$ in the fourth example 
$$ 
X_2=\{0,1,6,7,9,10,12,14,15,16,17,18,20,25,28,29,30,32,33,
37,39,43\}
$$
has cardinality 22. Consequently, its binary sequence is not balanced.
The correlation values of the sequence $X_2$, for the shifts 
in the range $0,1,\ldots,24$ are:
$$
49,1,1,1,1,1,1,-3,1,-3,1,1,-3,1,-3,-3,1,-3,1,1,-3,-3,1,1,-3.
$$
Thus the correlation values of $X_2$ occupy only three levels,
49, 1 and $-3$. Hence, this sequence is optimal but not balanced. The same is true for the fifth example. 

The block $X_2$ in the last example 
\begin{eqnarray*}
X_2&=&\{1,2,3,9,12,13,15,19,22,26,27,28,31,33,34,35,36,37,39,\\
&& 41,42,45,46,47,49,54,56,57,58,59\}
\end{eqnarray*}
has cardinality 30 and so its binary sequence $X_2$ is balanced.
The correlation values of the sequence $X_2$, for the shifts 
in the range $0,1,\ldots,30$ are:

\begin{eqnarray*}
&& 61,1,-3,-3,-3,-3,1,1,1,-3,1,-3,1,1,1,1,-3,1, \\
&& \quad  1,-3,-3,-3,-3,1,1,-3,-3,1,-3,-3,1.
\end{eqnarray*}

Hence, $X_2$ is a balanced optimal binary sequence of length 61.

\section{Propus difference families} \label{sec:PropDF}

In Table 3 we list propus difference families that we 
constructed by using the method of orbits. We only consider the cases where the subgroup $H$ is nontrivial. If each of the 
$k_i$ is the size of an $H$-invariant subset of $\bZ_v$, then 
we say that the parameter set is {\em H-feasible} (or just 
{\em feasible} when $H$ is known from the context). 
The case $v=67$ is omitted as it was treated 
separately in section \ref{sec:67}.

We can permute the $X_i$ and replace any $X_i$ with
its complement. When listing the propus difference families it 
is convenient to introduce some additional restrictions on 
the propus parameter sets (\ref{param-gen}). We shall assume that each $k_i\le v/2$, $k_2=k_3$ and that $k_1\ge k_4$.

In Table 3 below we first record the propus parameter set, and the subgroup $H$ of the multiplicative group of the finite field $\bZ_v$. Each of the three blocks $X_1$, $X_2=X_3$, $X_4$ is a union of orbits of $H$ acting on the additive group of $\bZ_v$. 
In order to specify which orbits constitute a block we
just list the representatives of these orbits.
As representative we choose the smallest integer in the orbit. 
For instance, $0$ is the unique representative of the trivial 
orbit $\{0\}$, and $1$ is the representative of the orbit $H$.

When two or more difference families are listed for the same 
parameter set, they are separated by a semicolon. 
When $k_1>k_4$ we have tried to find propus difference 
families with $X_1$ symmetric as well as those with $X_4$ symmetric. However, in some cases we did not succeed.

\begin{center} Table 3. Propus difference families with 
$v\equiv 1 \pmod{6}$ a prime.
\end{center}

\begin{verbatim}
(7;3,3,3,1;3), H={1,2,4}
[3],[3],[0]

(13;6,6,6,3;8), H={1,3,9}
[1,4],[4,7],[4]

(13;6,4,4,6;7), H={1,3,9}
[2,7],[0,4],[1,7]

(19;7,9,9,6;12), H={1,7,11}
[0,4,10],[2,4,5],[1,10]; [0,1,8],[1,4,10],[1,8]

(19;9,7,7,7;11), H={1,7,11}
[2,4,8],[0,5,10],[0,1,8]

(31;15,15,15,10;24), H={1,2,4,8,16}
[3,7,15],[1,3,15],[1,15]

(31;15,12,12,13;21), H={1,5,25}
[1,2,4,8,12],[2,4,8,11],[0,2,4,11,12]

(31;13,13,13,12;20), H={1,5,25}
[0,1,2,6,12],[0,2,6,8,11],[2,4,12,16];
[0,2,4,11,17],[0,3,8,11,17],[1,4,6,11]
 
(37;18,15,15,15;26), H={1,10,26}
[2,3,5,7,17,18],[1,3,7,17,21],[6,7,14,17,21]

(37;16,18,18,13;28), H={1,10,26}
[0,1,7,14,17,21],[1,2,6,9,14,21],[0,1,2,11,17]

(43;21,21,21,15;35), H={1,4,11,16,21,35,41}
[6,7,9],[1,6,9],[0,3,6]
 
(43;19,18,18,18;30), H={1,6,36}           
[0,2,4,9,14,19,20],[2,3,10,13,20,26],[3,4,10,13,20,21];
[0,5,7,9,10,20,21],[1,3,4,10,14,21],[1,3,5,7,13,21]

(43;18,21,21,16;33), H={1,6,36}
[1,5,7,10,13,26],[2,3,5,13,14,20,26],[0,1,7,9,19,20]

(49;22,24,24,18;39), H={1,18,30}
[0,1,6,7,8,9,13,16],[3,7,8,9,13,16,21,29],[3,6,8,12,16,29];
[0,2,7,8,13,16,19,26],[2,6,9,12,16,24,26,29],[1,3,7,8,19,21];
[0,1,3,4,12,13,16,24],[1,6,8,13,16,19,24,29],[1,4,6,16,19,26]

(49;22,22,22,19;36), H={1,18,30}
[0,4,6,7,9,13,19,26],[0,1,6,7,9,12,16,29],[0,1,6,7,16,19,21];
[0,3,4,6,7,12,19,29],[0,1,2,4,7,8,13,19],[0,1,3,7,8,19,21]

(61;30,26,26,26;47), H={1,9,20,34,58}
[2,6,8,10,23,26],[0,1,4,5,6,8],[0,3,5,6,10,12]

(61;30,25,25,30;49), H={1,9,20,34,58}
[4,5,10,12,13,26],[1,5,6,8,26],[2,4,10,12,13,26]

(61;30,25,25,30;49), H={1,13,47}                
[1,4,6,8,9,11,14,18,23,32],[0,6,7,8,14,22,23,27,28],
[1,3,4,6,7,8,9,11,28,36];

(61;25,30,30,25;49), H={1,13,47}
[0,3,6,7,8,18,22,23,31],[1,2,9,14,16,18,22,23,31,36],
[0,1,8,9,18,27,28,31,36]

The sequence of the second block below has only four correlation 
values, 61, 1, -3 and -11.
(61;28,28,28,24;47), H={1,13,47}
[0,1,3,4,14,16,18,23,31,32],[0,3,4,9,14,16,18,22,28,32],
[2,6,8,11,18,23,28,32]

(61;28,27,27,25;49), H={1,13,47}
[0,1,2,4,7,8,16,28,32,36],[1,2,7,8,9,12,16,27,36],
[0,2,7,12,16,27,28,31,36]

(73;36,36,36,28;63), H={1,8,64}
[3,5,6,11,12,21,25,26,27,33,35,43],
[3,4,9,14,17,18,21,26,34,35,42,43],[0,1,7,13,18,21,25,33,35,42]

(73;36,31,31,33;58), H={1,8,64}
[3,4,5,6,13,14,25,27,33,34,36,42],[0,2,3,5,9,18,21,26,27,35,42],
[1,5,7,11,18,21,27,33,34,42,43]

(73;31,36,36,30;60), H={1,8,64}
[0,2,7,11,12,13,17,18,26,35,42],
[3,5,6,12,14,18,21,26,27,33,34,35],[1,2,5,6,9,12,26,34,36,42]

(73;31,34,34,31;55), H={1,8,64}
[0,1,3,5,7,9,12,17,27,33,35],[0,1,2,5,9,11,12,18,21,27,36,43],
[0,1,3,9,18,21,26,27,35,36,42]

(73;31,36,36,30;60), H={1,8,64}
[0,1,4,14,17,21,26,34,36,42,43],
[2,3,4,7,12,14,25,27,35,36,42,43],[1,4,9,11,12,13,26,35,36,42]

(73;34,33,33,30;57), H={1,8,64}
[0,2,3,4,6,7,9,12,13,26,27,35],[1,2,5,6,7,12,17,21,25,26,35],
[2,4,6,7,11,17,18,25,26,36]

(157;78,78,78,66;143), 
H={1,14,16,39,46,67,75,93,99,101,108,130,153}
[2,3,7,9,11,13],[3,5,6,11,13,15],[0,3,4,5,7,13]

(307;153,153,153,136;288),
H={1,9,81,115,114,105,24,216,102,304,280,64,269,272,299,235,273}
[2,3,4,5,6,7,14,20,30],[4,5,7,12,14,28,30,31,49],
[2,6,7,10,21,28,30,31]
\end{verbatim}

The last two families have the same parameter sets as the 
corresponding Turyn propus families of the same lengths but 
they are not equivalent to them.

\section{Acknowledgements}
The research of the first author leading to these results has 
received funding from the Ministry of Education and Science of 
the Russian Federation according to the project part of the 
state funding assignment No 2.2200.2017/4.6.
The second author acknowledges generous support by NSERC. His 
work was enabled in part by support provided by the 
Shared Hierarchical Academic Research Computing Network (www.sharcnet.ca) and Compute/Calcul Canada (www.computecanada.ca).

\end{document}